\documentclass[11pt,reqno,table]{amsart}

\usepackage[T1]{fontenc}
\usepackage{amsmath}
\usepackage{amssymb}
\usepackage{amsthm}
\usepackage{amscd}
\usepackage{amsfonts}
\usepackage{mathrsfs}
\usepackage{mathtools}
\usepackage{stmaryrd}
\usepackage{algorithm, algorithmic}
\usepackage{wasysym}
\usepackage{caption}
\usepackage{subcaption}
\usepackage[dvipsnames]{xcolor}
\usepackage{microtype}

\usepackage{extarrows}
\usepackage{url}
\usepackage[colorlinks=true,hyperindex,linkcolor=magenta,pagebackref=false,citecolor=cyan]{hyperref}
\usepackage[alphabetic,lite,backrefs]{amsrefs}
\usepackage{tikz}
\usepackage{tikz-cd}
\usepackage{fullpage}
\usepackage[shortlabels]{enumitem}
\usepackage{standalone}

\DeclareMathAlphabet{\pazocal}{OMS}{zplm}{m}{n}

\newtheorem{lemma}{Lemma}[section]
\newtheorem{theorem}[lemma]{Theorem}
\newtheorem{prop}[lemma]{Proposition}
\newtheorem{cor}[lemma]{Corollary}

\newtheorem{defn}[lemma]{Definition}

\newtheorem*{theorem*}{Theorem}

\theoremstyle{remark}
\newtheorem{remark}[lemma]{Remark}
\newtheorem{example}[lemma]{Example}

\newcommand{\reg}{\operatorname{reg}}

\newcommand{\Spec}{\operatorname{Spec}}

\newcommand{\ann}{\operatorname{ann}}

\newcommand{\Pic}{\operatorname{Pic}}

\newcommand{\Nef}{\operatorname{Nef}}
\newcommand{\Eff}{\operatorname{Eff}}
\newcommand{\Fitt}{\operatorname{Fitt}}

\newcommand{\doot}{\bullet}

\renewcommand{\aa}{\mathbf a}
\newcommand{\bb}{\mathbf b}
\newcommand{\cc}{\mathbf c}
\newcommand{\CC}{\mathbf C}
\newcommand{\dd}{\mathbf d}

\newcommand{\ww}{\mathbf w}

\newcommand{\pp}{\mathbf p}
\newcommand{\PP}{\mathbf P}
\newcommand{\qq}{\mathbf q}

\newcommand{\cH}{\mathcal{H}}
\newcommand{\cI}{\mathcal{I}}

\renewcommand{\O}{\mathcal{O}}

\newcommand{\K}{\mathbb{K}}

\newcommand{\N}{\mathbb{N}}
\renewcommand{\P}{\mathbb{P}}

\newcommand{\Z}{\mathbb{Z}}

\newif\ifshownotes \shownotestrue 

\theoremstyle{theorem}
\newtheorem*{thm:reg-M-fg}{Theorem~\ref{thm:reg-M-fg}}
\newtheorem*{thm:powers}{Theorem~\ref{thm:powers}}

\title{Bounds on Multigraded Regularity}

\author{Juliette Bruce}
\address{Department of Mathematics, Dartmouth College, Hanover, NH}
\email{\href{mailto:juliette.bruce@dartmouth.edu}{juliette.bruce@dartmouth.edu}}

\author{Lauren Cranton Heller}
\address{Department of Mathematics, University of Nebraska, Lincoln, NE}
\email{\href{mailto:lheller2@unl.edu}{lheller2@unl.edu}}

\author{Mahrud Sayrafi}
\address{School of Mathematics, University of Minnesota, Minneapolis, MN}
\email{\href{mailto:mahrud@umn.edu}{mahrud@umn.edu}}

\subjclass[2020]{13D02, 14M25}

\begin{document}



\maketitle


\section{Introduction}

Castelnuovo--Mumford regularity is a measure of complexity for modules and sheaves in commutative algebra and algebraic geometry. It can be used as a bound on the degrees of syzygies of a module or the twists for which a sheaf is generated by global sections.

In \cite{maclaganSmith04}, Maclagan and Smith defined the multigraded Castelnuovo--Mumford regularity $\reg M\subseteq\Pic X$ of a $\Pic X$-graded module $M$ corresponding to a sheaf on a simplicial projective toric variety $X$\!.  We will specialize their definition to the case when $X$ is smooth.

The region $\reg M$ is invariant under translation by the nef cone of $X$, denoted by $\Nef X$. Extending intuition from the standard notion of regularity for $\PP^n$, one might expect that if $M$ is finitely generated then $\reg M$ could be specified by finitely many minimal elements with respect to $\Nef X$.  The following example shows this finite generation can fail even in simple situations.

\begin{example}\label{ex:not-fg}
	Let $M=S/\langle x_2,x_3\rangle$ be the coordinate ring of a single point on the Hirzebruch surface $\cH_t$ (see Example~\ref{ex:hirz}).  Since $\langle x_2,x_3\rangle$ is saturated $H_B^0(M)=0$.  Further, since the support of $\widetilde M$ has dimension $0$ we must have $H_B^i(M)=0$ for $i\geq 2$.  Thus $\reg M$ is determined entirely by $H_B^1(M)$, which vanishes exactly where the Hilbert function of $M$ agrees with its Hilbert polynomial.
	
	The Hilbert function of $M$ is equal to 1 inside the effective cone $\Eff\cH_t$ (see Figure~\ref{fig:hirzebruch}) and zero outside of it, so $\reg M = \Eff\cH_t$.  When $t>0$ this cone does not contain finitely many minimal elements with respect to $\Nef X$, as illustrated in Figure~\ref{fig:regularity-module-combined}.
  
\begin{figure}[h]
\newcommand{\makeaxes}{
  \path[use as bounding box] (-7.5,-1.5) rectangle (7.5,7.5);
  
  \fill[fill=blue!20] (0,0) -- (7,0) -- (7,7) -- (-7,7) -- (-7,7/2) -- cycle;
  \fill[fill=blue!35] (0,0) -- (7,0) -- (7,7) -- ( 0,7) -- cycle;}
\newcommand{\makegrid}{
  \draw[-, semithick] (-7,0)--(7,0);
  \draw[-, semithick] (0,-1)--(0,7);

  \foreach \x in {-7,...,7}
  \foreach \y in {-1,...,7}
    { \fill[gray,fill=gray] (\x,\y) circle (1.5pt); }}

\begin{tikzpicture}[scale=.33]
  \makeaxes
  \fill[fill=PineGreen!50] (0,0) -- (0,1) -- (-2,1) -- (-2,2) -- (-4,2) -- (-4,3) -- (-6,3) -- (-6,4) -- (-7,4) -- (-7,7) -- (7,7) -- (7,0) -- cycle;
  \makegrid
  \fill[TealBlue, fill=PineGreen] (0,0) circle (6pt);
  \draw[->, ultra thick, PineGreen] (0,0) -- (7,0);
  \draw[-, cap=round, ultra thick, PineGreen] (7,0) -- (0,0) -- (0,1) -- (-2,1) -- (-2,2) -- (-4,2) -- (-4,3) -- (-6,3) -- (-6,4) -- (-7,4);
\end{tikzpicture}
\hspace{0.5in}
\begin{tikzpicture}[scale=.33]
  \makeaxes
  \fill[fill=PineGreen!50] (0,7) -- (0,1) -- (7,1) -- (7,7) -- cycle;
  \makegrid
  \draw[->, ultra thick, PineGreen] (0,1) -- (0,7);
  \draw[->, ultra thick, PineGreen] (0,1) -- (7,1);
\end{tikzpicture}
\caption{
  On left, the multigraded regularity (green) of the module in Example~\ref{ex:not-fg} is an infinite staircase contained in a translate of the effective cone of $\cH_2$ (blue). \\
  On right, the multigraded regularity (green) of the ideal itself is contained in the nef cone of $\cH_2$ (dark blue).}
\label{fig:regularity-module-combined}
\end{figure}

\end{example}

Despite this apparent pathology we prove that when $M$ is faithful, meaning that $\ann M=0$, $\reg M$ is finitely generated for an arbitrary smooth projective toric variety $X$.

\begin{thm:reg-M-fg}
	Let $M$ be a finitely generated graded faithful $S$-module with $\widetilde M\neq 0$. Then $\reg M$ is contained in a translate of $\Nef X$. In particular, $\reg M$ has finitely many minimal elements.
\end{thm:reg-M-fg}

The proof produces an explicit bound by a translate of $\Nef X$, determined by the degrees of the generators of $M$ (see the figure in Example~\ref{ex:regularity-module}).  We use the idea that if the truncation $M_{\geq\dd}$ is not generated in a single degree $\dd$ then $M$ is not $\dd$-regular (see Theorem~\ref{thm:reg-S-fg} for a simpler case). 

Results from \cites{maclaganSmith04,botbolChardin17} have excluded particular degrees from the regularity region, implying that it cannot contain elements lying below these degrees in the partial order determined by $\Nef X$.  This does not preclude the existence of infinitely many incomparable elements.

Note that Theorem~\ref{thm:reg-M-fg} applies when $M=I_{Y}$ is the defining ideal of a subvariety $Y\subset X$ but not when $M=S/I_{Y}$ is the total coordinate ring of $Y$, and as in Example~\ref{ex:not-fg} we cannot expect $\reg(S/I_{Y})$ to be finitely generated in general. This is an issue which does not appear in the classical setting, where $\reg(I_Y)$ and $\reg(S/I_Y)$ differ by 1.  However, our result shows that in the primary case of historical interest---ideals and ideal sheaves of subvarieties---multigraded regularity is well-behaved.

It remains an interesting problem to characterize modules with regularity inside a translate of the nef cone.  We expect that the faithfulness hypothesis in Theorem~\ref{thm:reg-M-fg} can be weakened.   

As Example~\ref{ex:not-fg} shows, the regularity of an arbitrary finitely generated module may fail to be contained in all translates of $\Nef X$.  The regularity of the module in Example~\ref{ex:not-fg} is nevertheless contained in a translate of $\Eff X$.  We show in Proposition~\ref{prop:eff-bound} that this is true for all $M$. Hence the existence of a module whose regularity contains infinitely many minimal elements is a consequence of the difference between the effective and nef cones of $X$. This possibility highlights a theme from \cites{bruceCrantonhellerSayrafi21,berkeschKleinLoperYang22} that algebraic properties which coincide over projective spaces can diverge in higher Picard rank.

\subsection{Powers of Ideals}

Building on  the work of Swanson in \cite{swanson97}, Cutkosky--Herzog--Trung \cite{cutkoskyHerzogTrung99} and Kodiyalam \cite{kodiyalam00}  described the surprisingly predictable asymptotic behavior of Castelnuovo--Mumford regularity for powers of ideals on a projective space $\P^r$: given an ideal $I\subset \K[x_0,\ldots,x_r]$, there exist $d,e\in\Z$ such that for $n\gg0$ the regularity of $I^n$ satisfies \[ \reg\!\left(I^n\right) = dn+e. \]
Due to the importance of regularity as a measure of complexity for syzygies and its geometric interpretation in terms of the cohomology of coherent sheaves \cites{bertramEinLazarsfeld91,cutkoskyEinLazarsfeld01}, this phenomenon has received substantial attention \cites{geramitaGimiglianoPitteloud95,chandler97,smithSwanson97,romer01,trungWang05,bagheriChardinHa13}, focused mostly on projective spaces. See \cite{chardin13} for a survey.

As an application of our finite generation result we show how the asymptotic linearity of regularity for powers of ideals can be generalized to other smooth projective toric varieties. We bound the multigraded regularity region by describing regions ``inside'' and ``outside'' of $\reg\!\left(I^n\right)$ each of which translates linearly by a fixed vector as $n$ increases (see the figure in Example~\ref{ex:regularity-ideal}).  

The inner bound depends on the Betti numbers of the Rees ring $S[It]$, while the outer bound depends only on the degrees of the generators of $I$.  We use the partial order defined by $\Nef X$.

\begin{thm:powers}
  There exists a degree $\aa\in\Pic X$, depending only on $I$, such that for each integer $n>0$ and each pair of degrees $\qq_1,\qq_2\in\Pic X$ satisfying $\qq_1\geq\deg f_i\geq\qq_2$ for all generators $f_i$ of $I$, we have
	\[ n\qq_1+\aa+\reg S \subseteq \reg\!\left(I^n\right) \subseteq n\qq_2+\Nef X. \]
\end{thm:powers}

It is worth emphasizing that our result holds over smooth projective toric varieties with arbitrary Picard rank. Indeed, toric varieties of higher Picard rank introduce a wrinkle that is not present in existing asymptotic results on Castelnuovo--Mumford regularity: in general there are infinitely many possible regularity regions compatible with two given bounds.  (In contrast, when $\Pic X=\Z$, inner and outer bounds correspond to upper and lower bounds, respectively, with only finitely many integers between each pair.)

\subsection*{Outline}

The organization of the paper is as follows:
Section~\ref{sec:notation} introduces background results and our notation.
Section~\ref{sec:reg-S} shows that the multigraded regularity of $S$ lies inside $\Nef X$, in Theorem~\ref{thm:reg-S-fg}, and Section~\ref{sec:reg-M} shows that the multigraded regularity of a finitely generated faithful $S$-module is contained in an appropriate translate of $\Nef X$, in Theorem~\ref{thm:reg-M-fg}.  Section~\ref{sec:ideals} gives explicit inner and outer bounds for the multigraded regularity of powers of an ideal, in Theorem~\ref{thm:powers}.

\subsection*{Acknowledgments}

We thank Christine Berkesch and David Eisenbud for their helpful conversations and comments, and Daniel Erman for suggesting this problem, discussing it with us, and pointing out Example~\ref{ex:not-fg}. We also thank Dave Jensen, Tyler Kelly, and Hunter Spink for conversations which indirectly contributed to this project. The computer algebra system \texttt{Macaulay2}~\cite{M2}, in particular the package \texttt{NormalToricVarieties} by Gregory G.~Smith et~al., was indispensable in computing examples.

The first author is grateful for the support of the Mathematical Sciences Research Institute in Berkeley, California, where she was in residence for the 2020--2021 academic year. The first author was partially supported by the National Science Foundation under Award Nos. DMS-1440140, NSF FRG DMS-2053221, and NSF MSPRF DMS-2002239. The second author was partially supported by the NSF grant DMS-2001649 and the third author by the NSF grant DMS-2001101.

\section{Notation and Definitions}\label{sec:notation}

Throughout we work over a base field $\K$ and denote by $\N$ the set of non-negative integers. Let $X$ be a smooth projective toric variety determined by a fan.  The total coordinate ring of $X$ is a $\Pic(X)$-graded polynomial ring $S$ over $\K$ with an irrelevant ideal $B\subset S$.  Write $\Eff X$ for the monoid in $\Pic X$ generated by the degrees of the variables in $S$.

Fix minimal generators $\CC=(\cc_1,\ldots,\cc_r)$ for the monoid $\Nef X$ of classes in $\Pic X$ represented by numerically effective divisors.  For $\lambda\in\Z^r$, write $\lambda\cdot\CC$ to represent the linear combination $\lambda_1\cc_1+\cdots+\lambda_r\cc_r\in\Pic X$, and similarly for other tuples in $\Pic X$.  Write $|\lambda|$ for the sum $\lambda_1+\cdots+\lambda_r$.  We use a partial order on $\Pic X$ induced by $\Nef X$: given $\aa,\bb\in\Pic X$, we write $\aa\leq\bb$ when $\bb-\aa\in\Nef X$.

\begin{example}\label{ex:hirz}
  The Hirzebruch surface $\cH_t = \P\left(\O_{\P^{1}}\oplus\O_{\P^{1}}(t)\right)$ is a smooth projective toric variety whose associated fan, shown left in Figure~\ref{fig:hirzebruch}, has rays $(1,0)$, $(0,1)$, $(-1,t)$, and $(0,-1)$. For each ray there is a corresponding prime torus-invariant divisor.  In particular, the total coordinate ring of $\cH_t$ is the polynomial ring $S=\K[x_0,x_1,x_2,x_3]$ and its irrelevant ideal is $B = \langle x_0,x_2\rangle\cap\langle x_1,x_3\rangle$.
  
\begin{figure}[h]
  \begin{tikzpicture}[scale=.5]
    \path[use as bounding box] (-3,-3) rectangle (3,3);

    \fill[fill=blue!10] (0,0) -- (3, 0) -- ( 3,  3) -- ( 0, 3) -- cycle;
    \fill[fill=blue!25] (0,0) -- (0, 3) -- (-3/2,3) -- ( 0, 0);
    \fill[fill=blue!40] (0,0) -- (0,-3) -- (-3, -3) -- (-3, 3) -- (-3/2,3) -- (0,0);
    \fill[fill=blue!55] (0,0) -- (3, 0) -- ( 3, -3) -- ( 0,-3) -- cycle;

    \draw[line width=1pt,black] (0,0) -- (-3/2,3);
    \draw[line width=1pt,black] (0,0) -- ( 0, -3);
    \draw[line width=1pt,black] (0,0) -- ( 3,  0);
    \draw[line width=1pt,black] (0,0) -- ( 0,  3);

    \draw[line width=1.5pt,black,-stealth] (0,0) -- ( 1, 0) node[anchor=north west]{$\rho_0$}; 
    \draw[line width=1.5pt,black,-stealth] (0,0) -- ( 0,-1) node[anchor=south east]{$\rho_3$}; 
    \draw[line width=1.5pt,black,-stealth] (0,0) -- (-1, 2) node[anchor=north east]{$\rho_2$}; 
    \draw[line width=1.5pt,black,-stealth] (0,0) -- ( 0, 1) node[anchor=south west]{$\rho_1$}; 
  \end{tikzpicture}
  \hspace{1in}
  \begin{tikzpicture}[scale=.6]
    \path[use as bounding box] (-3,-2) rectangle (3,3);

    \fill[fill=blue!20] (0,0) -- (3,0) -- (3,3) -- (-3,3) -- (-3,3/2) -- cycle;
    \fill[fill=blue!35] (0,0) -- (3,0) -- (3,3) -- (0,3) -- cycle;

    \draw [thin, gray,] (-3,0) -- (3,0); 
    \draw [thin, gray,] (0,-1) -- (0,3); 

    \foreach \x in {-3,...,3}{
      \foreach \y in {-1,...,3}{
        \node[circle,draw=gray, fill=gray,inner sep=0.7pt] at (\x,\y){ };
    }}

    \draw[line width=1.5pt,black,-stealth] (0,0) -- ( 1,0) node[anchor=south west]{$x_0,x_2$};
    \draw[line width=1.5pt,black,-stealth] (0,0) -- (-2,1) node[anchor=south west]{$x_1$};
    \draw[line width=1.5pt,black,-stealth] (0,0) -- ( 0,1) node[anchor=south west]{$x_3$};
  \end{tikzpicture}
  \caption{Left: fan of $\cH_2$. Right: the cones $\Nef\cH_2$ (dark blue) and $\Eff\cH_2$ (blue).}\label{fig:hirzebruch}
\end{figure}

  Choosing a basis for $\Pic\cH_t \cong \Z^2$, the grading on $S$ can be given as $\deg x_0 = \deg x_2 = (1,0)$, $\deg x_1 = (-t,1)$, and $\deg x_3=(0,1)$.  The effective and nef cones are illustrated on the right.
\end{example}

For a $\Pic(X)$-graded $S$-module $M$ and $\dd\in\Pic X$, denote by $M_{\geq\dd}$ the submodule of $M$ generated by all elements of degrees $\dd'$ satisfying $\dd'\geq\dd$ (c.f.~\cite[Def.~5.1]{maclaganSmith04}).  Denote by $\widetilde M$ the quasi-coherent sheaf on $X$ associated to $M$, as in \cite[\S3]{cox95}.

We now recall the notion of multigraded Castelnuovo--Mumford regularity introduced by Maclagan and Smith.

\begin{defn}[{c.f.~\cite[Def.~1.1]{maclaganSmith04}}]
	Let $M$ be a graded $S$-module. For $\dd\in\Pic X$, we say $M$ is \emph{$\dd$-regular} if the following hold:
  \begin{enumerate}
	\item  $H^i_B(M)_\bb=0$ for all $i>0$ and all $\bb\in\bigcup_{|\lambda|=i-1}\left(\dd-\lambda\cdot\CC+\Nef X\right)$ where $\lambda\in\N^r$.
	\item  $H^0_B(M)_\bb=0$ for all $\bb\in\bigcup_j \left(\dd+\cc_j+\Nef X\right)$.
  \end{enumerate}
	We write $\reg M$ for the set of $\dd$ such that $M$ is $\dd$-regular.
\end{defn}

The following example shows that the degree of a minimal generator of an ideal does not bound its multigraded regularity on an arbitrary toric variety.

\begin{example}\label{ex:degs-dont-bound}
	Let $I=\langle x_0x_3,x_0x_2,x_1x_2\rangle$ be an ideal in the total coordinate ring of the Hirzebruch surface $\cH_t$, with notation as in Example~\ref{ex:hirz}. A local cohomology computation verifies that $I$ is $(1,1)$-regular. However $x_0x_2$ is a minimal generator with $\deg(x_0x_2)=(2,0) \not\leq (1,1)$.
\end{example}

The existence of a similar example with $H_B^0(M)\neq 0$ was noted by Maclagan and Smith, who asked whether $B$-torsion was necessary in \cite[\S5]{maclaganSmith04}. Example~\ref{ex:degs-dont-bound} shows that it is not.

\section{Regularity of the Total Coordinate Ring}\label{sec:reg-S}

We begin with the case $M = S$, showing that the pathology seen in Example~\ref{ex:not-fg}---a regularity region contained in no translate of $\Nef X$---does not occur for the total coordinate ring of a smooth projective toric variety. In particular we show that $\reg S\subseteq\Nef X$.

In \cite[Prob.~6.12]{maclaganSmith04}, Maclagan and Smith asked for a combinatorial characterization of toric varieties $X$ such that $\Nef X\subseteq\reg S$.  Theorem~\ref{thm:reg-S-fg} below shows that when $X$ is smooth and projective, $\Nef X\subseteq\reg S$ is in fact equivalent to the a priori stronger condition that $\reg S = \Nef X$. It still remains an interesting question to characterize such toric varieties. For instance, the only Hirzebruch surface with this property is $\cH_1$.

\begin{theorem}\label{thm:reg-S-fg}
	Using the notation from Section~\ref{sec:notation}, we have $\reg S \subseteq \Nef X$. In particular, $\reg S$ contains finitely many minimal elements.
\end{theorem}
\begin{proof}
	Take $\dd\in\reg S$.  By \cite[Thm.~5.4]{maclaganSmith04} the truncation $S_{\geq\dd}$ is generated by the monomials of $S_\dd$, so there is a surjection $S_\dd\otimes_\K S\to S_{\geq\dd}(\dd)$ which sheafifies to a surjection $S_\dd\otimes\O\to\O(\dd)$. Hence $\O(\dd)$ is generated by global sections, so by \cite[Thm.~6.3.11]{coxLittleSchenck11} $\dd$ is nef.

	An application of Dickson's lemma (e.g.~\cite[\S2.4 Thm.~5]{coxLittleOShea15}), suggested by Will Sawin~\cite{sawinMO}, shows that $\reg S$ has finitely many minimal elements, finishing the proof.

\begin{lemma}\label{lem:mins-in-cone}
	A subset $V\subseteq\Nef X$ contains finitely many minimal elements with respect to $\leq$ on $\Pic X$.
\end{lemma}

 Elements of $V$ can be written as linear combinations $\lambda\cdot\CC$ of the monoid generators of $\Nef X$.  The minimal elements of $V$ must have coefficients $\lambda\in\N^r$ that are minimal in the component-wise partial order on $\N^r$.  By Dickson's lemma only finitely many possible coefficients exist.
\end{proof}

\begin{example}
	The multigraded regularity of the total coordinate ring of the Hirzebruch surface $\cH_2$ is contained in the nef cone of $\cH_2$, as illustrated in Figure~\ref{fig:regularity-hirzebruch}.
	
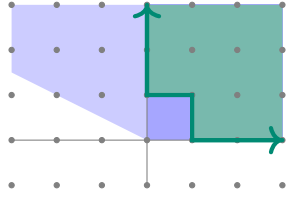
\begin{figure}[h]
  \begin{tikzpicture}[scale=.6]
    \path[use as bounding box] (-3,-1) rectangle (3,3);

    \fill[fill=blue!20] (0,0) -- (3,0) -- (3,3) -- (-3,3) -- (-3,3/2) -- cycle;
    \fill[fill=blue!35] (0,0) -- (3,0) -- (3,3) -- (0,3) -- cycle;
    \fill[fill=PineGreen!50] (1,0) -- (3,0) -- (3,3) -- (0,3) -- (0,1) -- (1,1) -- cycle;

    \draw [thin, gray,] (-3,0) -- (3,0); 
    \draw [thin, gray,] (0,-1) -- (0,3); 

    \foreach \x in {-3,...,3}{
      \foreach \y in {-1,...,3}{
        \node[circle,draw=gray, fill=gray,inner sep=0.7pt] at (\x,\y){ };
    }}

    \draw[->, ultra thick, PineGreen] (1,0) -- (3,0);
    \draw[->, ultra thick, PineGreen] (0,1) -- (0,3);
    \draw[-, cap=round, ultra thick, PineGreen] (3,0) -- (1,0) -- (1,1) -- (0,1) -- (0,3);
  \end{tikzpicture}
  \caption{The regularity of $S$ (dark green) is contained in $\Nef\cH_2$ (dark blue).}\label{fig:regularity-hirzebruch}
\end{figure}

\end{example}

Though we do not directly use Theorem~\ref{thm:reg-S-fg} in the next section, we do rely on the idea of the proof.  For an arbitrary module $M$, if $\dd\in\reg M$ then the truncation $M_{\geq\dd}$ is generated in a single degree $\dd$, meaning that $\widetilde M(\dd)$ is globally generated.  This no longer immediately implies that $\dd$ is nef, but Lemma~\ref{lem:extra-monomials} below connects the difference between $\dd$ and the degrees of the generators of $M$ to monomials in truncations of $S$ itself.

We also use the chamber complex of the rays of $\Eff X$, which is described in \cite[\S2]{maclaganSmith04}.  By definition, this chamber complex is the coarsest fan with support $\Eff X$ which refines all triangulations of the degrees of the variables of $S$.  It partitions $\Eff X$ into cones that govern many geometric properties of $\Spec S$, including its GIT quotients, birational geometry, and Hilbert polynomials (c.f.\ \cite[Ch.~14-15]{coxLittleSchenck11}, \cite[\S5]{heringKuronyaPayne06}).

For our purposes we need only the existence of a strongly convex rational polyhedral fan that covers $\Eff X$ and contains $\Nef X$ as a cone.  We will refer to the maximal cones as chambers and the codimension one cones as walls. In particular, $\Nef X$ is a chamber.

\begin{lemma}\label{lem:extra-monomials}
	Let $\Gamma$ be a chamber of $\Eff X$ other than $\Nef X$, and let $\aa_1,\ldots,\aa_n\in\Pic X$.  If $\aa_i\in\Gamma\setminus\Nef X$ for all $i$, then there exist monomials $m_i\in S_{\geq\aa_i}$ such that $\prod_i m_i$ is not generated by the monomials of $S_{\sum\aa_i}$.
\end{lemma}
\begin{proof}
	Since $\Gamma$ and $\Nef X$ intersect at most in a wall of $\Gamma$ and no $\aa_i$ lies in $\Gamma\cap\Nef X$, their sum $\bb = \sum\aa_i$ must also be in $\Gamma\setminus\Nef X$.  Consider the multiplication maps
	\[\begin{tikzcd}
	{S_\bb\otimes_\K S} & {S(\bb)} \\
	& {\bigotimes_\K S_{\geq\aa_i}(\aa_i)}.
	\arrow["\varphi", from=1-1, to=1-2]
	\arrow["\psi"', from=2-2, to=1-2]
	\end{tikzcd}\]
	Suppose the proposition is false.  Then the image of $\psi$ must be contained in the image of $\varphi$, else we could choose $(m_i)\in\bigotimes_\K S_{\geq\aa_i}(\aa_i)$ with image not generated by the monomials of $S_\bb$. Note that each $S_{\geq\aa_i}(\aa_i)$ sheafifies to $\O(\aa_i)$, so sheafifying the entire diagram gives
	\[\begin{tikzcd}
	{S_\bb\otimes \O} & {\O(\bb)} \\
	& {\O(\bb)}.
	\arrow["\varphi", from=1-1, to=1-2]
	\arrow["\psi"', from=2-2, to=1-2]
	\end{tikzcd}\]
	In particular, the image of $\psi$ is still contained in the image of $\varphi$.  Since $\psi$ sheafifies to an isomorphism, $\varphi$ sheafifies to a surjection. This implies $\bb\in\Nef X$, which is a contradiction.
\end{proof}


\section{Regularity of Faithful Modules}\label{sec:reg-M}

The goal of this section is to prove that the multigraded regularity of a faithful module has only finitely many minimal elements.

Proposition~\ref{prop:eff-bound} shows that the regularity of an arbitrary finitely generated module is contained in some translate of $\Eff X$.  Under the additional assumption that $M$ is faithful, i.e.\ that $\ann M=0$, Proposition~\ref{prop:nef-bound} shows we can also eliminate degrees that are in a translate of $\Eff X$ but not $\Nef X$.

\begin{prop}\label{prop:eff-bound}
	Let $M$ be a finitely generated graded $S$-module with $\widetilde M\neq 0$. Suppose the degrees of all minimal generators of $M$ are contained in $\Eff X$. Then $\reg M \subseteq \Eff X$.
\end{prop}
\begin{proof}
	Take $\dd\in\reg M$ and suppose for contradiction that $\dd\not\in\Eff X$. The degree $\dd$ part $M_\dd$ generates $M_{\geq\dd}$ by \cite[Thm.~5.4]{maclaganSmith04}. By hypothesis all elements of $M$ have degrees inside $\Eff X$, so $M_\dd = 0$ and thus $M_{\geq\dd}=0$. The modules $M$ and $M_{\geq\dd}$ define the same sheaf by \cite[Lem.~6.8]{maclaganSmith04}, so $M_{\geq\dd}=0$ contradicts $\widetilde M\neq 0$.
\end{proof}

\begin{prop}\label{prop:nef-bound}
	Let $M$ be a finitely generated graded faithful $S$-module with $\widetilde M\neq 0$.  Suppose $\Gamma$ is a chamber of $\Eff X\setminus\Nef X$. If $\dd-\deg f_i\in\Gamma\setminus\Nef X$ for all generators $f_i$ of $M$, then $M$ is not $\dd$-regular.
\end{prop}
\begin{proof}
	Assume on the contrary that $M$ is $\dd$-regular.  Let $\aa_i=\dd-\deg f_i$ for each $i$.  By choice of $\dd$ we have $\aa_i\in\Gamma\setminus\Nef X$. Hence by Lemma~\ref{lem:extra-monomials} there exist monomials $m_i\in S_{\geq\aa_i}$ such that $\prod_i m_i$ is not generated by the monomials of $S_{\sum\aa_i}$. Consider the elements $m_i f_i \in M_{\geq\dd}$.

	Since $M$ is $\dd$-regular, the degree $\dd$ part $M_\dd$ generates $M_{\geq\dd}$ by \cite[Thm.~5.4]{maclaganSmith04}. Let $g_1,\dots,g_s$ with $\deg g_j = \dd$ be generators for $M_{\geq\dd}$. Thus we must have relations
	\[ m_i f_i = \sum_j b_{i,j} g_j = \sum_j b_{i,j} \left(\sum_k a_{j,k} f_k\right) = \sum_k c_{i,k} f_k \]
	for some $b_{i,j},a_{j,k},c_{i,k}\in S$ with $\deg b_{i,j} = \deg m_i - \aa_i$ and $\deg a_{j,k} = \aa_k$. These relations form a partial presentation matrix
  \begin{align}\label{eq:A}
	  A =
    \begin{bmatrix}
	    m_1 & 0 & \cdots & 0 \\
	    0 & m_2 & \cdots & 0 \\
	    \vdots  & \vdots &
      \ddots  & \vdots \\
	    0 & 0 & \cdots & m_n
	  \end{bmatrix} -
	  \begin{bmatrix}
	    c_{1,1} & c_{2,1} & \cdots & c_{n,1} \\
	    c_{1,2} & c_{2,2} & \cdots & c_{n,2} \\
	    \vdots  & \vdots  & \ddots & \vdots  \\
	    c_{1,n} & c_{2,n} & \cdots & c_{n,n}
	  \end{bmatrix}.
  \end{align}
	for $M$.  In particular, $\det(A)\in\Fitt_0 M\subseteq\ann M$ by \cite[Prop.~20.7]{eisenbud95}, so $\det(A) M = 0$.

	Since $\ann M=0$ we must have $\det(A) = 0$, but this is impossible: note that $\det(A)$ contains the monomial $m = \prod_i m_i$ and that $\det(A) \in m + I$ for $I = \prod_k\langle c_{1,k},c_{2,k},\dots,c_{n,k}\rangle$, then observe that $I \subseteq \prod_k\langle a_{1,k},a_{2,k},\dots,a_{n,k}\rangle\subseteq S\otimes_\K S_{\sum\aa_k}$ since $\deg a_{j,k} = \aa_k$. Hence $\det(A) = 0$ implies $m \in I\subseteq S\otimes_\K S_{\sum\aa_k}$ and contradicts our choice of $m_i$.
\end{proof}

\begin{remark}
	Example~\ref{ex:not-fg} shows that there are non-faithful modules which do not satisfy the conclusion of Theorem~\ref{thm:reg-M-fg}.  In practice, however, we only need some choice of $m_i$ as in Lemma~\ref{lem:extra-monomials} which lead to a contradiction.  Given a specific toric variety, it may be possible to verify directly that the element $\det A$ from \eqref{eq:A} cannot annihilate $M$ in some cases where $M$ is not faithful.
\end{remark}

We will use the following technical lemma about the walls of $\Nef X$ to find a vector satisfying the hypotheses of Proposition~\ref{prop:nef-bound}.

\begin{lemma}\label{lem:chambers}
	Given $\aa_1,\ldots,\aa_n\in\Nef X$ and $\dd\in\Eff X\setminus\Nef X$, there exists a chamber $\Gamma$ sharing a wall $W$ with $\Nef X$ and $\ww$ in the relative interior of $W$ such that $\dd+\ww\in\Gamma$ and $\dd+\ww\in\aa_i+\Gamma$ for all $i$.
\end{lemma}
\begin{proof}
  Consider the cone $P$ defined by all rays of $\Nef X$ in addition to a primitive element along $\dd$.  Since $\Nef X\subsetneq P$, at least one wall $W$ of $\Nef X$ must be in the interior of $P\subseteq\Eff X$.  Let $\Gamma$ be the chamber across $W$ from $\Nef X$.  Since $\dd\notin\Nef X$, for each $\ww\in W$ we have $\dd+\ww\notin\Nef X$.
  
\usetikzlibrary{patterns}
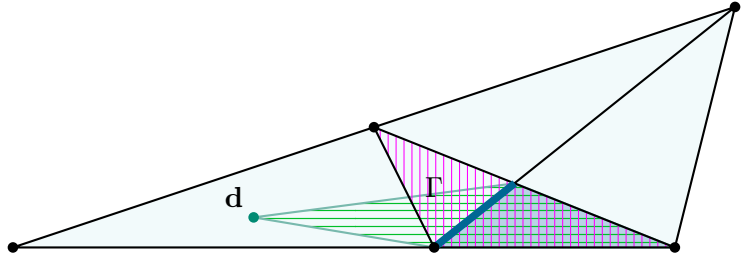
\begin{figure}[h]
\newcommand{\makeaxes}{
  \path[use as bounding box] (-7,0) rectangle (7,7);
  \draw[join=round, thick, black] (-12,0) -- (10,0) -- (0,4) -- (2,0) -- (12,8) -- cycle;
  \draw[join=round, thick, black] (10,0) -- (12,8) -- cycle;
}
\begin{tikzpicture}[scale=.4]
  \fill[fill=Turquoise!5] (-12,0) -- (10,0) -- (12,8) -- cycle;
  
  \fill[fill=MidnightBlue!20] (10,0) -- (421/90,192/90) -- (2,0) -- cycle;
  
  \fill[pattern=horizontal lines, pattern color=PineGreen] (10,0) -- (421/90,192/90) -- (-4,1) -- (2,0) -- cycle;
  \draw[join=round, thick, PineGreen!50] (10,0) -- (421/90,192/90) -- (-4,1) -- (2,0) -- cycle;
  
  \fill[pattern=vertical lines, pattern color=Magenta] (10,0) -- (0,4) -- (2,0) -- cycle;
  \makeaxes
  \draw[join=round, line width=1mm, MidnightBlue] (421/90,192/90) -- (2,0);
  \fill[fill=black] (-12,0) circle (5pt);
  \fill[fill=black] (  2,0) circle (5pt);
  \fill[fill=black] ( 10,0) circle (5pt);
  \fill[fill=black] ( 12,8) circle (5pt);
  \fill[fill=black] (  0,4) circle (5pt);
  \fill[fill=PineGreen] (-4,1) circle (5pt) node[anchor=south east]{$\dd$};
  \fill[fill=PineGreen] (2,2) node{$\Gamma$};
\end{tikzpicture}
\caption{A section of a hypothetical chamber complex with $P$ (green, horizontal) and $Q$ (red, vertical) inside $\Eff X$. The chamber $\Nef X$ and its wall $W$ are in blue.}\label{fig:chambers}
\end{figure}

  Now consider the cone $Q$ defined by all supporting hyperplanes of $\Nef X$ and $\Gamma$ except the hyperplane containing $W$.  Since $W$ is in the intersection of the open half-spaces defining $Q$, it lies in the interior of $Q$.  Therefore we can find $\ww$ in the relative interior of $W\subset Q$ so that $\dd+\ww\in\aa_i+Q\subseteq\aa_i+(\Gamma\cup\Nef X)$ for all $i$.  By hypothesis $\aa_i+\Nef X\subseteq\Nef X$ so $\dd+\ww\notin\aa_i+\Nef X$.  Hence $\dd+\ww\in\aa_i+\Gamma$ for all $i$.
\end{proof}

\begin{theorem}\label{thm:reg-M-fg}
	Let $M$ be a finitely generated graded faithful $S$-module with $\widetilde M\neq 0$. Suppose the degrees of all minimal generators of $M$ are contained in $\Nef X$. Then $\reg M \subseteq \Nef X$. In particular, $\reg M$ has finitely many minimal elements.
\end{theorem}
\begin{proof}
	Suppose there exists $\dd\in\reg M\setminus\Nef X$. Since $M$ satisfies the hypothesis of Proposition~\ref{prop:eff-bound}, we can assume that $\dd\in\Eff X$.  Using Lemma~\ref{lem:chambers}, we can find $\ww$ in the relative interior of a wall separating $\Nef X$ and an adjacent chamber $\Gamma$ such that $\dd+\ww\in\Gamma$ and $\dd+\ww\in\deg f_i+\Gamma$ for all $i$.  It follows from Proposition~\ref{prop:nef-bound} that $\dd+\ww\notin\reg M$, which is a contradiction because $\ww\in\Nef X$ and $\reg M$ is invariant under positive translation by $\Nef X$.

	The conclusion that $\reg M$ has finitely many minimal elements follows from Lemma~\ref{lem:mins-in-cone}.
\end{proof}

\begin{cor}\label{cor:outer-bound}
	Let $M$ be a finitely-generated faithful $S$-module.  If $\deg f_i\in\bb+\Nef X$ for all generators $f_i$ of $M$ then $\reg M\subseteq\bb+\Nef X$.
\end{cor}

\begin{example}\label{ex:regularity-module}
  Consider the Hirzebruch surface $\cH_2$, with notation from Example~\ref{ex:hirz}, and let $M$ be the torsion-free module with presentation
  \[\begin{tikzcd}[column sep=8em, ampersand replacement=\&]
	  {S(3, -3) \oplus S(2, -2) \oplus S(1, -2)} \& {S(0, -4).}
	  \arrow["\;\;\;\;{\begin{pmatrix} x_0^5 x_1 & x_1^2 x_2^6 & x_1^2 x_2^5 \end{pmatrix}^T}"', from=1-2, to=1-1]
  \end{tikzcd}\]
  Since the degrees of the generators are contained in $(-3,2) + \Nef\cH_2$, by Corollary~\ref{cor:outer-bound} the multigraded regularity of $M$ is contained in a translate of the nef cone, illustrated in Figure~\ref{fig:regularity-module}.
  
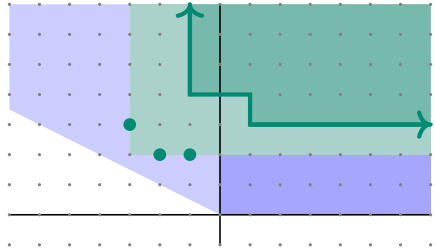
\begin{figure}[h]
\newcommand{\makeaxes}{
  \path[use as bounding box] (-7.5,-1.5) rectangle (7.5,7.5);
  
  \fill[fill=blue!20] (0,0) -- (7,0) -- (7,7) -- (-7,7) -- (-7,7/2) -- cycle;
  \fill[fill=blue!35] (0,0) -- (7,0) -- (7,7) -- ( 0,7) -- cycle;}
\newcommand{\makegrid}{
  \draw[-, semithick] (-7,0)--(7,0);
  \draw[-, semithick] (0,-1)--(0,7);

  \foreach \x in {-7,...,7}
  \foreach \y in {-1,...,7}
    { \fill[gray,fill=gray] (\x,\y) circle (1.5pt); }}

\begin{tikzpicture}[scale=.33]
  \makeaxes
  \fill[fill=PineGreen!30] (-3,7) -- (-3,2) -- (7,2) -- (7,7) -- cycle;
  \fill[fill=PineGreen!50] (-1,7) -- (-1,4) -- (1,4) -- (1,3) -- (7,3) -- (7,7) -- cycle;
  \makegrid
  \fill[fill=PineGreen] (-3,3) circle (6pt);
  \fill[fill=PineGreen] (-2,2) circle (6pt);
  \fill[fill=PineGreen] (-1,2) circle (6pt);
  \draw[->, ultra thick, PineGreen] (-1,4) -- (-1,7);
  \draw[->, ultra thick, PineGreen] ( 1,3) -- ( 7,3);
  \draw[-, cap=round, ultra thick, PineGreen] (-1,7) -- (-1,4) -- (1,4) -- (1,3) -- (7,3);
\end{tikzpicture}
\caption{The multigraded regularity (dark green) of the module $M$ is contained in a translate $(-3,2) + \Nef X$ (light green) of the nef cone of $\cH_2$ (dark blue).}\label{fig:regularity-module}
\end{figure}

\end{example}

\section{Powers of Ideals and Multigraded Regularity}\label{sec:ideals}

Throughout this section let $I=\langle f_1,\ldots,f_s\rangle \subseteq S$ be an ideal and let $\PP$ be the vector with coordinates $\pp_i=\deg f_i\in\Pic X$.  We are interested in the asymptotic behavior of the multigraded regularity of $I^n$ as $n$ increases.  In particular, we prove the following theorem:

\begin{theorem}\label{thm:powers}
  There exists a degree $\aa\in\Pic X$, depending only on $I$, such that for each integer $n>0$ and each pair of degrees $\qq_1,\qq_2\in\Pic X$ satisfying $\qq_1\geq\pp_i\geq\qq_2$ for all $i$, we have
	\[ n\qq_1+\aa+\reg S \subseteq \reg\!\left(I^n\right) \subseteq n\qq_2+\Nef X. \]
\end{theorem}
\begin{proof}
	The inner bound will follow from Proposition~\ref{prop:inner-bound}.  The outer bound follows from Corollary~\ref{cor:outer-bound} by noting that $\deg\prod_{j=1}^nf_{i_j}=\sum_{j=1}^n\pp_{i_j}\in n\qq_2+\Nef X$ for all products of $n$ choices of generators of $I$, and such products generate $I^n$.
\end{proof}

\begin{example}\label{ex:regularity-ideal}
  Let $I = \langle x_0x_3, x_1^2x_2^4\rangle$ and $J = \langle x_3, x_0^3x_1\rangle$ be two ideals in the total coordinate ring of the Hirzebruch surface $\cH_2$, with notation as in Example~\ref{ex:hirz}. Figure~\ref{fig:regularity-ideal} shows the multigraded regularity of powers of $I$ and $J$ along with the bounds from Theorem~\ref{thm:powers}.
  
\begin{figure}[ht]
\newcommand{\makegrid}{
  \path[use as bounding box] (-1.5,-1.5) rectangle (11.5,11.5);

  \draw[-, semithick] (-1,0)--(11,0);
  \draw[-, semithick] (0,-1)--(0,11);

  \foreach \x in {-1,...,11}
  \foreach \y in {-1,...,11}
    { \fill[gray,fill=gray] (\x,\y) circle (1.5pt); }}

\begin{tikzpicture}[scale=.25]
  \fill[fill=PineGreen!30] (0,11) -- (0,1) -- (11,1) -- (11,11) -- cycle;
  \fill[fill=PineGreen!50] (1,11) -- (1,3) -- (2,3) -- (2,2) -- (11,2) -- (11,11) -- cycle;
  \fill[fill=PineGreen!70] (2,11) -- (2,4) -- (3,4) -- (3,3) -- (11,3) -- (11,11) -- cycle;
  \makegrid
  \fill[fill=PineGreen] (1,1) circle (6pt);
  \fill[fill=PineGreen] (0,2) circle (6pt);
  \draw[->, ultra thick, PineGreen] (1,3) -- (1,11);
  \draw[->, ultra thick, PineGreen] (2,2) -- (11,2);
  \draw[-, cap=round, ultra thick, PineGreen] (1,11) -- (1,3) -- (2,3) -- (2,2) -- (11,2);
  \node at (4.5,-2) {{\footnotesize $\reg\!\left(I\right)$}};
\end{tikzpicture}
\quad
\begin{tikzpicture}[scale=.25]
  \fill[fill=PineGreen!30] (0,11) -- (0,2) -- (11,2) -- (11,11) -- cycle;
  \fill[fill=PineGreen!50] (2,11) -- (2,4) -- (11,4) -- (11,11) -- cycle;
  \fill[fill=PineGreen!70] (3,11) -- (3,6) -- (4,6) -- (4,5) -- (11,5) -- (11,11) -- cycle;
  \makegrid
  \fill[fill=PineGreen] (2,2) circle (6pt);
  \fill[fill=PineGreen] (1,3) circle (6pt);
  \fill[fill=PineGreen] (0,4) circle (6pt);
  \draw[->, ultra thick, PineGreen] (2,4) -- (2,11);
  \draw[->, ultra thick, PineGreen] (2,4) -- (11,4);
  \draw[-, cap=round, ultra thick, PineGreen] (2,11) -- (2,4) -- (11,4);
  \node at (4.5,-2) {{\footnotesize $\reg\!\left(I^2\right)$}};
\end{tikzpicture}
\quad
\begin{tikzpicture}[scale=.25]
  \fill[fill=PineGreen!30] (0,11) -- (0,3) -- (11,3) -- (11,11) -- cycle;
  \fill[fill=PineGreen!50] (3,11) -- (3,6) -- (11,6) -- (11,11) -- cycle;
  \fill[fill=PineGreen!70] (4,11) -- (4,8) -- (5,8) -- (5,7) -- (11,7) -- (11,11) -- cycle;
  \makegrid
  \fill[fill=PineGreen] (3,3) circle (6pt);
  \fill[fill=PineGreen] (2,4) circle (6pt);
  \fill[fill=PineGreen] (1,5) circle (6pt);
  \fill[fill=PineGreen] (0,6) circle (6pt);
  \draw[->, ultra thick, PineGreen] (3,6) -- (3,11);
  \draw[->, ultra thick, PineGreen] (3,6) -- (11,6);
  \draw[-, cap=round, ultra thick, PineGreen] (3,11) -- (3,6) -- (11,6);
  \node at (4.5,-2) {{\footnotesize $\reg\!\left(I^3\right)$}};
\end{tikzpicture}
\quad
\begin{tikzpicture}[scale=.25]
  \fill[fill=PineGreen!30] (0,11) -- (0,4) -- (11,4) -- (11,11) -- cycle;
  \fill[fill=PineGreen!50] (4,11) -- (4,8) -- (11,8) -- (11,11) -- cycle;
  \fill[fill=PineGreen!70] (5,11) -- (5,10) -- (6,10) -- (6,9) -- (11,9) -- (11,11) -- cycle;
  \makegrid
  \fill[fill=PineGreen] (4,4) circle (6pt);
  \fill[fill=PineGreen] (3,5) circle (6pt);
  \fill[fill=PineGreen] (2,6) circle (6pt);
  \fill[fill=PineGreen] (1,7) circle (6pt);
  \fill[fill=PineGreen] (0,8) circle (6pt);
  \draw[->, ultra thick, PineGreen] (4,8) -- (4,11);
  \draw[->, ultra thick, PineGreen] (4,8) -- (11,8);
  \draw[-, cap=round, ultra thick, PineGreen] (4,11) -- (4,8) -- (11,8);
  \node at (4.5,-2) {{\footnotesize $\reg\!\left(I^4\right)$}};
\end{tikzpicture}

\vspace{1em}
\begin{tikzpicture}[scale=.25]
  \fill[fill=PineGreen!30] (0,11) -- (0,1) -- (11,1) -- (11,11) -- cycle;
  \fill[fill=PineGreen!50] (0,11) -- (0,2) -- (2,2) -- (2,1) -- (11,1) -- (11,11) -- cycle;
  \fill[fill=PineGreen!70] (2,11) -- (2,3) -- (3,3) -- (3,2) -- (11,2) -- (11,11) -- cycle;
  \makegrid
  \fill[fill=PineGreen] (0,1) circle (6pt);
  \fill[fill=PineGreen] (1,1) circle (6pt);
  \draw[->, ultra thick, PineGreen] (0,2) -- (0,11);
  \draw[->, ultra thick, PineGreen] (2,1) -- (11,1);
  \draw[-, cap=round, ultra thick, PineGreen] (0,11) -- (0,2) -- (2,2) -- (2,1) -- (11,1);
  \node at (4.5,-2) {{\footnotesize $\reg\!\left(J\right)$}};
\end{tikzpicture}
\quad
\begin{tikzpicture}[scale=.25]
  \fill[fill=PineGreen!30] (0,11) -- (0,2) -- (11,2) -- (11,11) -- cycle;
  \fill[fill=PineGreen!50] (0,11) -- (0,4) -- (1,4) -- (1,3) -- (3,3) -- (3,2) -- (11,2) -- (11,11) -- cycle;
  \fill[fill=PineGreen!70] (3,11) -- (3,4) -- (4,4) -- (4,3) -- (11,3) -- (11,11) -- cycle;
  \makegrid
  \fill[fill=PineGreen] (0,2) circle (6pt);
  \fill[fill=PineGreen] (1,2) circle (6pt);
  \fill[fill=PineGreen] (2,2) circle (6pt);
  \draw[->, ultra thick, PineGreen] (0,4) -- (0,11);
  \draw[->, ultra thick, PineGreen] (3,2) -- (11,2);
  \draw[-, cap=round, ultra thick, PineGreen] (0,11) -- (0,4) -- (1,4) -- (1,3) -- (3,3) -- (3,2) -- (11,2);
  \node at (4.5,-2) {{\footnotesize $\reg\!\left(J^2\right)$}};
\end{tikzpicture}
\quad
\begin{tikzpicture}[scale=.25]
  \fill[fill=PineGreen!30] (0,11) -- (0,3) -- (11,3) -- (11,11) -- cycle;
  \fill[fill=PineGreen!50] (0,11) -- (0,5) -- (2,5) -- (2,4) -- (4,4) -- (4,3) -- (11,3) -- (11,11) -- cycle;
  \fill[fill=PineGreen!70] (4,11) -- (4,5) -- (5,5) -- (5,4) -- (11,4) -- (11,11) -- cycle;
  \makegrid
  \fill[fill=PineGreen] (0,3) circle (6pt);
  \fill[fill=PineGreen] (1,3) circle (6pt);
  \fill[fill=PineGreen] (2,3) circle (6pt);
  \fill[fill=PineGreen] (3,3) circle (6pt);
  \draw[->, ultra thick, PineGreen] (0,5) -- (0,11);
  \draw[->, ultra thick, PineGreen] (4,3) -- (11,3);
  \draw[-, cap=round, ultra thick, PineGreen] (0,11) -- (0,5) -- (2,5) -- (2,4) -- (4,4) -- (4,3) -- (11,3);
  \node at (4.5,-2) {{\footnotesize $\reg\!\left(J^3\right)$}};
\end{tikzpicture}
\quad
\begin{tikzpicture}[scale=.25]
  \fill[fill=PineGreen!30] (0,11) -- (0,4) -- (11,4) -- (11,11) -- cycle;
  \fill[fill=PineGreen!50] (0,11) -- (0,7) -- (1,7) -- (1,6) -- (3,6) -- (3,5) -- (5,5) -- (5,4) -- (11,4) -- (11,11) -- cycle;
  \fill[fill=PineGreen!70] (5,11) -- (5,6) -- (6,6) -- (6,5) -- (11,5) -- (11,11) -- cycle;
  \makegrid
  \fill[fill=PineGreen] (0,4) circle (6pt);
  \fill[fill=PineGreen] (1,4) circle (6pt);
  \fill[fill=PineGreen] (2,4) circle (6pt);
  \fill[fill=PineGreen] (3,4) circle (6pt);
  \fill[fill=PineGreen] (4,4) circle (6pt);
  \draw[->, ultra thick, PineGreen] (0,7) -- (0,11);
  \draw[->, ultra thick, PineGreen] (5,4) -- (11,4);
  \draw[-, cap=round, ultra thick, PineGreen] (0,11) -- (0,7) -- (1,7) -- (1,6) -- (3,6) -- (3,5) -- (5,5) -- (5,4) -- (11,4);
  \node at (4.5,-2) {{\footnotesize $\reg\!\left(J^4\right)$}};
\end{tikzpicture}
\caption{The inner (dark green) and outer (light green) bounds for powers of $I$ and $J$. The circles correspond to the degrees of the generators of each power.}\label{fig:regularity-ideal}
\end{figure}
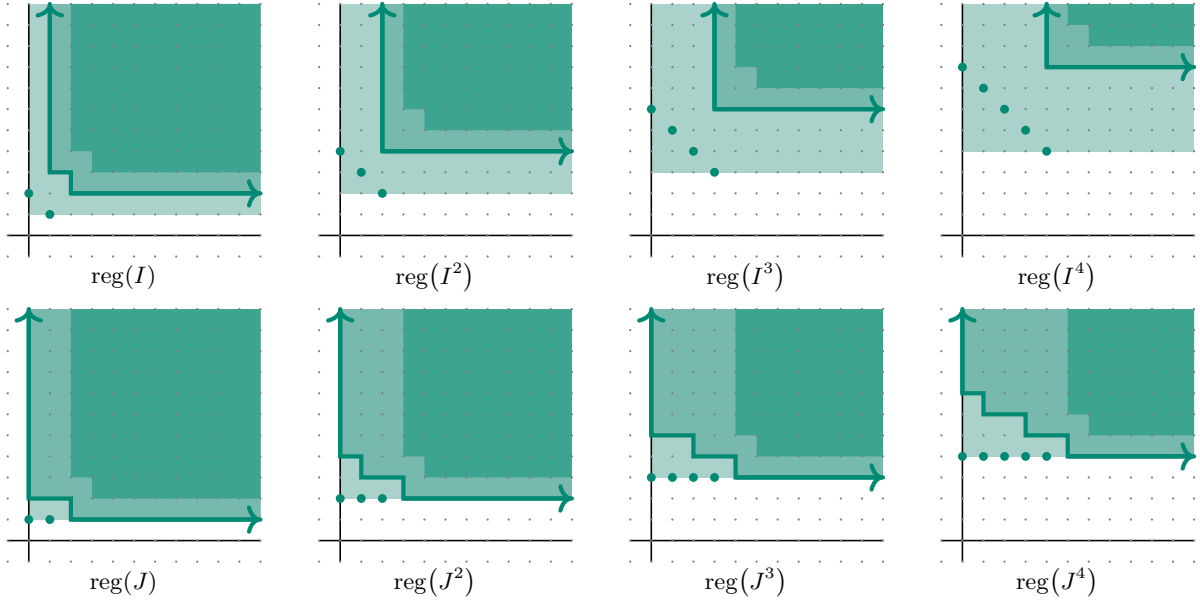

\end{example}

\begin{remark}
	If $\qq_2$ is not nef, then the bounds in Theorem~\ref{thm:powers} will not increase with $n$ in the partial order on $\Pic X$.  We can see that this behavior is necessary by taking $I$ to be a principal ideal generated outside of $\Nef X$.
\end{remark}

\subsection{The Rees Ring}\label{sec:rees}

One way to find a subset of the regularity of a module is by using its multigraded Betti numbers.  In order to describe $\reg\!\left(I^n\right)$, we would thus like a uniform description of the Betti numbers of $I^n$ for all $n$.  For this purpose, consider the multigraded Rees ring of $I$:
\[ S[It] \coloneqq  \bigoplus_{n\geq0} I^nt^n \subseteq S[t], \]
which is a $\Pic(X)\times\Z$-graded noetherian ring with $\deg ft^k=(\deg f,k)$ for $f\in S$. Let $R=S[T_1,\ldots,T_s]$ be the $\Pic(X)\times\Z$-graded ring with $\deg(T_i)=(\deg f_i,1)=(\pp_i,1)$. Notice that there is a surjective map of graded $S$-algebras:
\[\begin{tikzcd}[row sep = -0.25em]
R \rar[two heads]{} & S[It] \\
T_i \rar[mapsto]{} & f_i\, t
\end{tikzcd}\]
Since $R$ is a finitely generated standard graded algebra over $S$, taking a single degree of a finitely generated $R$-module in the auxiliary $\Z$ grading yields a finitely generated $S$-module.

\begin{defn}
	For a $\Pic(X)\times\Z$-graded $R$-module $M$, define $M^{(n)}$ to be the $\Pic(X)$-graded $S$-module
	\[ M^{(n)}\coloneqq \bigoplus_{\aa\in\Pic X} M_{(\aa,n)}. \]
\end{defn}

Following \cite{kodiyalam00}, we record three important properties of this operation.

\begin{lemma}\label{lem:properties-of-(n)}
  Consider the functor $-^{(n)}\colon M\mapsto M^{(n)}$ from the category of $\Pic(X)\times\Z$-graded $R$-modules to the category of $\Pic(X)$-graded $S$-modules.
\begin{enumerate}[(i)]
	\item\label{it:exact} $-^{(n)}$ is an exact functor.
	\item\label{it:powers} $S[It]^{(n)}\cong I^n$.
	\item\label{it:R-to-S} $R(-\aa,-b)^{(n)} \cong R^{(n-b)}(-\aa) \cong \bigoplus_{|\nu|=n-b} S(-\nu\cdot\PP-\aa)$ where $\nu\in\N^s$.
\end{enumerate}
\end{lemma}

Since $S[It]$ is a finitely generated module over the polynomial ring $R$, it has a finite free resolution. Applying $-^{(n)}$ gives a resolution by \ref{it:exact}, which has cokernel $I^n$ by \ref{it:powers} and whose terms are finitely generated free $S$-modules by \ref{it:R-to-S}.  Thus we can constrain the Betti numbers of $I^n$ in terms of those of $S[It]$.

\subsection{Regularity of Powers of Ideals}\label{sec:reg-In}

Given a description of the Betti numbers of $I^n$ in terms of $n$, we obtain an inner bound on $\reg\!\left(I^n\right)$ using the following lemma.

\begin{lemma}\label{lem:MS-variant}
	If $F_\bullet$ is a finite free resolution for $M$ with $F_j=\bigoplus_i S(-\aa_{i,j})$ and $H_B^0(M)=0$ then
	\begin{align}\label{eq:inclusion}
		\bigcap_{i,j}\bigcup_{|\lambda|=j} (\aa_{i,j}-\lambda\cdot\CC+\reg S)\subseteq\reg M
	\end{align}
	where $\CC=(\cc_1,\ldots,\cc_r)$ is the sequence of nef generators for $X$ and the union is over $\lambda\in\N^r$.
\end{lemma}

\begin{remark}
	This result amounts to switching the union and intersection in the statement of \cite[Cor.~7.3]{maclaganSmith04} for modules with $H_B^0(M)=0$, which increases the size of the subset by allowing a different choice of $\lambda$ for each $i,j$.
\end{remark}
\begin{proof}
	Fix $\dd$ in the left hand side of \eqref{eq:inclusion} and consider the hypercohomology spectral sequence for $F_\bullet$ (see \cite[Thm.~4.14]{bruceCrantonhellerSayrafi21} for a description of this spectral sequence).  We must show that $M$ is $\dd$-regular, meaning that $H^k_B(M)_{\dd-\mu\cdot\CC}=0$ for all $k$ and all $\mu$ with $|\mu|=k-1$.  Since $F_\bullet$ is a resolution for $M$, a diagonal of our spectral sequence converges to $H^k_B(M)$.  Thus it is sufficient to prove that this entire diagonal vanishes in degree $\dd-\mu\cdot\CC$, i.e.\ that
	\begin{align}\label{eq:summands}
		H^{k+j}_B(F_j)_{\dd-\mu\cdot\CC} = \bigoplus_i H^{k+j}_B(S(-\aa_{i,j}))_{\dd-\mu\cdot\CC} = 0
	\end{align}
	for all $j$.  This is satisfied for $k=0$ by hypothesis.  Now fix $k>0$, $\mu$, $j$, and $i$.  By choice of $\dd$ we have $\dd\in \aa_{i,j}-\lambda\cdot\CC+\reg S$ for some $\lambda$ with $|\lambda|=j$, so that $\dd-\aa_{i,j}+\lambda\cdot\CC\in\reg S$.  Call this degree $\dd'$, and let $\cc'=(\lambda+\mu)\cdot\CC$, where $|\lambda+\mu|=k+j-1$.  Then by the definition of the regularity of $S$ we have $H^{k+j}_B(S)_{\dd'-\cc'}=0$
	where
	\[ \dd'-\cc' = \dd-\aa_{i,j}+\lambda\cdot\CC-(\lambda+\mu)\cdot\CC = \dd-\mu\cdot\CC.\]
	Hence each summand in \eqref{eq:summands} is zero for $k>0$, as desired.
\end{proof}

\begin{prop}\label{prop:inner-bound}
	There exists a degree $\aa\in\Pic X$, depending only on the Rees ring of $I$, such that for each integer $n>0$ and degree $\qq\in\Pic X$ satisfying $\qq\geq\deg f_i$ for all homogeneous generators $f_i$ of $I$, we have
	\[ n\qq+\aa+\reg S \subseteq \reg\!\left(I^n\right). \]
\end{prop}
\begin{proof}
	\newcommand{\capsub}{\substack{i,j \\ |\nu| = n-b_{i,j}}}
	Let $F_{\doot}$ be a minimal $\Pic(X)\times\Z$-graded free resolution of $S[It]$ as an $R$-module, and write $F_j=\bigoplus_i R(-\aa_{i,j},-b_{i,j})$ for $\aa_{i,j}\in\Pic X$ and $b_{i,j}\in\Z$.
	By Lemma~\ref{lem:properties-of-(n)}, applying the $-^{(n)}$ functor to $F_{\doot}$ yields a (potentially non-minimal) resolution of $S[It]^{(n)}\cong I^n$ consisting of free $S$-modules
	\[ F_j^{(n)}
	\cong \bigoplus_i R(-\aa_{i,j}, -b_{i,j})^{(n)}
	\cong \bigoplus_i \left[ \bigoplus_{|\nu|=n-b_{i,j}}S(-\nu\cdot\PP-\aa_{i,j}) \right], \]
	where $\PP=(\deg f_1,\dots,\deg f_s)$ is the sequence of degrees of the homogeneous generators $f_i$ of $I$. From this Lemma~\ref{lem:MS-variant} gives the following bound on the regularity of $I^n$:
	\begin{align}\label{eq:proof-1}
		\bigcap_{\capsub} \bigcup_{|\lambda|=j} \left[ \nu\cdot\PP + \aa_{i,j} - \lambda\cdot\CC+\reg S \right]
		\subseteq \reg\!\left(I^n\right).
	\end{align}
	Note that $b_{0,0}=0$, as $S[It]$ is a quotient of $R$, and thus $b_{i,j}\geq 0$ for all $i,j$, as $R$ is positively graded in the $\Z$ coordinate.
	
	Take $\aa\in\Pic X$ so that $\aa\geq\aa_{i,j}$ for all $i,j$.  There are only finitely many $\aa_{i,j}$ because $S[It]$ is a finitely generated $R$-module and $R$ is noetherian. We may now simplify the left hand side of \eqref{eq:proof-1} by noting three things: (i) for all $|\lambda|=j$ and all $j$ we have $\reg S\subseteq-\lambda\cdot\CC+\reg S$, (ii) if $|\nu|=n-b_{i,j}$ then $(n-b_{i,j})\qq\in \nu\cdot\PP+\reg S$, and (iii) for all $i$ and all $j$ we have $n\qq+\aa\in(n-b_{i,j})\qq+\aa_{i,j}+\reg S$. 
	
	Combining these facts gives that
	\begin{align*}
		\reg\!\left(I^n\right)
		&\supseteq \bigcap_{\capsub} \bigcup_{|\lambda|=j} \left[ \nu\cdot\PP + \aa_{i,j} - \lambda\cdot\CC + \reg S \right] \\
		&\supseteq \bigcap_{\capsub} \left[ \nu\cdot\PP + \aa_{i,j} + \reg S \right] \\
		&\supseteq \hspace{1.25em} \bigcap_{i,j} \left[ (n-b_{i,j})\qq + \aa_{i,j} + \reg S \right] \\
		&\supseteq \quad n\qq+\aa+\reg S.
	\end{align*}
\end{proof}

A related problem is characterizing the asymptotic behavior of regularity for symbolic powers of $I$. Note that the symbolic Rees ring of $I$ is not necessarily noetherian (see \cite{grifoSeceleanu21}, for instance), so our argument for the existence of the degree $\aa$ in the proof of Proposition~\ref{prop:inner-bound} does not work in this case. More generally, if $\cI = \{I_n\}$ is a filtration of ideals, then one may ask for sufficient conditions so that $\reg\!\left(I_n\right)$ is uniformly bounded.



\begin{bibdiv}
\begin{biblist}

\bib{bagheriChardinHa13}{article}{
  author={Bagheri, Amir},
  author={Chardin, Marc},
  author={H\`a, Huy T\`ai},
  title={The eventual shape of Betti tables of powers of ideals},
  journal={Math. Res. Lett.},
  volume={20},
  date={2013},
  number={6},
  pages={1033--1046},
}

\bib{berkeschKleinLoperYang22}{article}{
  author={Berkesch, Christine},
  author={Klein, Patricia},
  author={Loper, Michael C.},
  author={Yang, Jay},
  title={Homological and combinatorial aspects of virtually Cohen-Macaulay
    sheaves},
  journal={Trans. London Math. Soc.},
  volume={9},
  date={2022},
  number={1},
  pages={413--434},
}

\bib{bertramEinLazarsfeld91}{article}{
  author={Bertram, Aaron},
  author={Ein, Lawrence},
  author={Lazarsfeld, Robert},
  title={Vanishing theorems, a theorem of {Severi}, and the equations defining projective varieties},
  journal={J. Amer. Math. Soc.},
  volume={4},
  date={1991},
  number={3},
  pages={587--602},
}

\bib{botbolChardin17}{article}{
	author={Botbol, Nicol\'as},
	author={Chardin, Marc},
	title={Castelnuovo Mumford regularity with respect to multigraded ideals},
	journal={J. Algebra},
	volume={474},
	date={2017},
	pages={361--392},
	issn={0021-8693},
	review={\MR{3595796}},
	doi={10.1016/j.jalgebra.2016.11.017},
}

\bib{bruceCrantonhellerSayrafi21}{article}{
  author = {Bruce, Juliette},
  author = {Cranton Heller, Lauren},
  author = {Sayrafi, Mahrud},
  title  = {Characterizing Multigraded Regularity on Products of Projective Spaces},
  date   = {2021},
  note   = {ArXiv pre-print: \url{https://arxiv.org/abs/2110.10705}}
}

\bib{chandler97}{article}{
  author={Chandler, Karen A.},
  title={Regularity of the powers of an ideal},
  journal={Comm. Algebra},
  volume={25},
  date={1997},
  number={12},
  pages={3773--3776},
}

\bib{chardin13}{article}{
  author={Chardin, Marc},
  title={Powers of ideals: Betti numbers, cohomology and regularity},
  conference={title={Commutative algebra},},
  book={publisher={Springer},},
  date={2013},
  pages={317--333},
}

\bib{cox95}{article}{
  author  = {Cox, David A.},
  title   = {The homogeneous coordinate ring of a toric variety},
  journal = {J. Algebraic Geom.},
  volume  = {4},
  date    = {1995},
  number  = {1},
  pages   = {17--50},
}

\bib{coxLittleOShea15}{book}{
  author={Cox, David A.},
  author={Little, John B.},
  author={O'Shea, Donal},
  title={Ideals, varieties, and algorithms},
  series={Undergraduate Texts in Mathematics},
  edition={4},
  publisher={Springer},
  date={2015},
  pages={xvi+646},
}

\bib{coxLittleSchenck11}{book}{
  author={Cox, David A.},
  author={Little, John B.},
  author={Schenck, Henry K.},
  title={Toric varieties},
  series={Graduate Studies in Mathematics},
  volume={124},
  publisher={AMS, Providence, RI},
  date={2011},
  pages={xxiv+841},
}

\bib{cutkoskyEinLazarsfeld01}{article}{
  author={Cutkosky, Steven Dale},
  author={Ein, Lawrence},
  author={Lazarsfeld, Robert},
  title={Positivity and complexity of ideal sheaves},
  journal={Math. Ann.},
  volume={321},
  date={2001},
  number={2},
  pages={213--234},
}

\bib{cutkoskyHerzogTrung99}{article}{
  author={Cutkosky, S. Dale},
  author={Herzog, J\"{u}rgen},
  author={Trung, Ng\^{o} Vi\^{e}t},
  title={Asymptotic behaviour of the Castelnuovo-Mumford regularity},
  journal={Compositio Mathematica},
  volume={118},
  date={1999},
  number={3},
  pages={243--261},
}

\bib{eisenbud95}{book}{
  author={Eisenbud, David},
  title={Commutative algebra},
  series={Graduate Texts in Mathematics},
  volume={150},
  publisher={Springer},
  date={1995},
  pages={xvi+785},
}

\bib{geramitaGimiglianoPitteloud95}{article}{
  author={Geramita, Anthony V.},
  author={Gimigliano, Alessandro},
  author={Pitteloud, Yves},
  title={Graded Betti numbers of some embedded rational $n$-folds},
  journal={Math. Ann.},
  volume={301},
  date={1995},
  number={2},
  pages={363--380},
}

\bib{M2}{misc}{
  author={Grayson, Daniel~R.},
  author={Stillman, Michael~E.},
  title={Macaulay2, a software system for research in algebraic geometry},
  note={Available at \url{https://faculty.math.illinois.edu/Macaulay2/}},
  label={M2},
}

\bib{grifoSeceleanu21}{article}{
  author={Grifo, Elo\'{\i}sa},
  author={Seceleanu, Alexandra},
  title={Symbolic Rees algebras},
  conference={
    title={Commutative algebra},
  },
  book={
    publisher={Springer},
  },
  date={2021},
  pages={343--371},
}

\bib{heringKuronyaPayne06}{article}{
  author={Hering, Milena},
  author={K\"{u}ronya, Alex},
  author={Payne, Sam},
  title={Asymptotic cohomological functions of toric divisors},
  journal={Adv. Math.},
  volume={207},
  date={2006},
  number={2},
  pages={634--645},
}

\bib{kodiyalam00}{article}{
  author={Kodiyalam, Vijay},
  title={Asymptotic behaviour of Castelnuovo-Mumford regularity},
  journal={Proc. Amer. Math. Soc.},
  volume={128},
  date={2000},
  number={2},
  pages={407--411},
}

\bib{maclaganSmith04}{article}{
  author={Maclagan, Diane},
  author={Smith, Gregory G.},
  title={Multigraded Castelnuovo-Mumford regularity},
  journal={J. Reine Angew. Math.},
  volume={571},
  date={2004},
  pages={179--212},
}


\bib{romer01}{article}{
  author={R\"{o}mer, Tim},
  title={Homological properties of bigraded algebras},
  journal={Illinois J. Math.},
  volume={45},
  date={2001},
  number={4},
  pages={1361--1376},
}

\bib{sawinMO}{misc}{
  title  = {Generalization of Dickson's Lemma},
  author = {Sawin, Will},
  note   = {URL: https://mathoverflow.net/q/383015 (version: 2021-02-03)},
  eprint = {https://mathoverflow.net/q/383015},
  organization={MathOverflow}
}

\bib{smithSwanson97}{article}{
  author={Smith, Karen E.},
  author={Swanson, Irena},
  title={Linear bounds on growth of associated primes},
  journal={Comm. Algebra},
  volume={25},
  date={1997},
  number={10},
  pages={3071--3079},
}

\bib{swanson97}{article}{
   author={Swanson, Irena},
   title={Powers of ideals: primary decompositions, Artin-Rees lemma and regularity},
   journal={Math. Ann.},
   volume={307},
   date={1997},
   number={2},
   pages={299--313},
}

\bib{trungWang05}{article}{
  author={Trung, Ng\^{o} Vi\^{e}t},
  author={Wang, Hsin-Ju},
  title={On the asymptotic linearity of Castelnuovo-Mumford regularity},
  journal={J. Pure Appl. Algebra},
  volume={201},
  date={2005},
  number={1-3},
  pages={42--48},
}

\end{biblist}
\end{bibdiv}
\end{document}